\documentstyle{amsart}
\input{bull-art}
\newcommand{\glue}{\tau}
\newcommand{\skin}{\sigma}

\newlength{\figboxwidth}             
\newcommand{\arrow}{\rightarrow}
\newcommand{\bdry}{\partial}

\newcommand{\compos}{\circ}

\newcommand{\dirsum}{\oplus}

\newcommand{\includesin}{\hookrightarrow}
\newcommand{\isom}{\cong}

\newcommand{\EQ}{\;\;=\;\;}

\newcommand{\st}{\;\: : \;\:}         

\newcommand{\zbar}{{\overline{z}}}

\newcommand{\cx}{{\Bbb C}}
\newcommand{\half}{{\Bbb H}}
\newcommand{\integers}{{\Bbb Z}}

\newcommand{\reals}{{\Bbb R}}

\newcommand{\chat}{\widehat{\cx}}

\newcommand{\zed}{\integers}

\newcommand{\Teic}{{\mathrm{Teich}}}

\newtheorem{theorem}{Theorem}[section]

\newcommand{\cM}{{\cal M}}

\title{ Riemann surfaces and the geometrization\\
 of 3-manifolds}

\author{	Curt McMullen }
\date{October 28, 1991. Presented to
		the Symposium on Frontiers of Mathematics,
		New York, 14--19 December 1988}
\thanks{
		Research partially supported by an 
		NSF Postdoctoral Fellowship}
\address{Department of Mathematics, Princeton University, 
Princeton,
New Jersey  08544}
\curraddr{Department of Mathematics, University of 
California, Berkeley, 
California  94720
  	}
\email{ctm@@math.berkeley.edu}
\begin{document}
\def\currentvolume{27}
\def\currentissue{2}
\def\currentyear{1992}
\def\currentmonth{October}
\def\copyrightyear{1992}
\def\currentpages{207-216}
\maketitle

\section{Introduction}

About a decade ago Thurston proved 
that a vast collection of 3-manifolds carry metrics of 
constant
negative curvature.  
These manifolds are thus elements of {\em hyperbolic 
geometry}, 
as natural as Euclid's regular polyhedra.
For a closed manifold, Mostow rigidity assures that a 
hyperbolic
structure is unique when it exists \cite{Mostow:book}, so 
topology
and geometry mesh harmoniously in dimension 3.

This remarkable theorem applies to all 3-manifolds, which 
can
be built up in an inductive way from 3-balls, i.e.,
{\em Haken} manifolds.  Thurston's construction of a 
hyperbolic
structure is also inductive.  At the inductive step
one must find the right geometry on an open 3-manifold
so that its ends may be glued together.  
Using quasiconformal deformations, the gluing problem
can be formulated as a fixed-point problem for
a map of Teichm\"uller space to itself.
Thurston proposes to find the fixed point by
{\em iterating} this map.

Here we outline Thurston's construction and sketch
a new proof that the 
iteration converges.  Our argument rests on a
result entirely in the theory of Riemann surfaces:
an extremal quasiconformal mapping can
be relaxed (isotoped to a map of lesser dilatation)
when lifted to a sufficiently large covering space
(e.g., the universal cover).  This contraction gives
an immediate estimate for the contraction of Thurston's
iteration.  

A detailed account of these results appears in
\cite{McMullen:amenable, McMullen:iter};
therefore, we have adopted a more informal tone here.
An expository sequel relating
hyperbolic manifolds and iterated rational maps
appears in \cite{McMullen:icm90}.

The setup for Thurston's theorem is discussed at
length in \cite{Morgan:uniformization}, and\break
Thurston's proof
is contained in the series [Th3, Th5, Th6, Th7], etc.
See also \cite{Thurston:survey}
for a survey and an introduction to the more general
geometrization conjecture and \cite{Thurston:book} for
more on hyperbolic geometry and Kleinian groups, as well as
the existence of hyperbolic structures on many manifolds 
which are not Haken.

Background in Teichm\"uller theory can be found, e.g., in
\cite{Bers:London:survey,
Bers:ts:survey, Gardiner:book}.

\section{Teichm\"uller theory and quasiconformal maps}
\label{sec:qc}

In this section we state results from 
\cite{McMullen:amenable} concerning contraction of 
canonical mappings between Teichm\"uller spaces.

Let $f : X \arrow Y$ be an orientation-preserving
homeomorphism between two Riemann surfaces $X$ and $Y$.

If $f$ is smooth then we can measure its conformal
distortion as follows.
First, consider a real-linear isomorphism $L : \cx \arrow 
\cx$.
Think of the complex structure on $\cx$ as being recorded
by the family of circles centered at the origin, since these
are the orbits of multiplication by the unit complex 
numbers. 
The image of these circles under $L$ is a family of ellipses
of constant eccentricity.  The ratio of major to minor axes
gives a number $K(L) \ge 1$ called the {\em dilatation} of 
$L$.

We define the dilatation $K(f)$ as the supremum of
the dilatation of the derivative $Df : T_x \arrow 
T_{f(x)}$, 
over all tangent spaces $T_x$ to $X$.

The dilatation of $f$ is $1$ if and only if 
$f$ is conformal.
If $K(f) < \infty$, $f$ is {\em quasiconformal}.

{\em Technical remark.} 
The natural degree of smoothness to
require is that the distributional derivatives of $f$ lie
in $L^2$. 

Assuming $f$ is quasiconformal, we can try to adjust it by
isotopy to minimize its dilatation.  By compactness of maps
with bounded dilatation, at least one
{\em extremal} (dilatation minimizing) quasiconformal map 
exists in each isotopy class.  
For example, if $X$ and $Y$ are tori, any extremal
map is an {\em affine stretch} (it takes the form
$z \mapsto A z + B \zbar + C$ in 
the universal covers of $X$ and $Y$).  

{\em Teichm\"uller's theorem} describes the extremal maps
when $X$ and $Y$ are hyperbolic Riemann surfaces of finite
area (equivalently, surfaces of negative Euler 
characteristic
obtained from compact surfaces by possibly removing a finite
number of points).  In each isotopy class there is a 
{\em unique} extremal {\em Teichm\"uller map}. 
Away from a finite number of singularities,
the map is again an affine stretch in appropriate local
coordinates.
\vskip 8pt

\noindent{\bf Relaxation of quasiconformal maps.}
Now consider quasiconformal mappings on
the unit disk $\Delta \EQ \{z \st |z|<1\}$. 

Any quasiconformal map $f : \Delta \arrow \Delta$
extends to a homeomorphism of $S^1$, the boundary of the 
disk.
Of course $f$ is always isotopic to the identity; to
obtain an interesting extremal problem, we require
that the isotopy {\em fix} the values of $f$ on $S^1$.

Our first result states that a Teichm\"uller mapping is
no longer extremal when lifted to the universal cover.

\begin{theorem}
\label{thm:not_extremal}
Let $f : X \arrow Y$ be a
Teichm\"uller mapping between hyperbolic Riemann
surfaces of finite type.  Then the map
$\tilde{f} : \Delta \arrow \Delta$ obtained by lifting $f$
to the universal covers of $X$ and $Y$
is not extremal among quasiconformal maps with
the same boundary values \RM(unless $f$ is conformal\RM).
\end{theorem}

This result can be reformulated in terms of natural maps
between {\em Teichm\"uller spaces}.  Given a
Riemann surface $X$, the Teichm\"uller space Teich($X$)
consists of equivalence classes of data $(f : X \arrow 
X_1)$ 
where $X_1$ is another Riemann surface and $f$ 
is a quasiconformal map.  Given two points $(f : X \arrow 
X_1)$
and $(g : X \arrow X_2)$, the {\em Teichm\"uller distance}
$d(X_1,X_2)$ is defined to be $\inf \log(K(h))$ where $h$
ranges over all quasiconformal maps $h : X_1 \arrow X_2$, 
which
are isotopic to $g\compos f^{-1}$ (rel ideal boundary).
Thus $d(X_1,X_2) = 0$ exactly when there is a conformal
isomorphism in the correct isotopy class, and in this
case we consider $(f,X_1)$ and $(g,X_2)$ to represent the 
same point
in Teichm\"uller space.

Now given a covering $Y \arrow X$, there is a natural 
inclusion
$\Teic(X) \includesin \Teic(Y)$; one simply forms the 
corresponding
covering space of each Riemann surface quasiconformally
equivalent to $X$.  For the special case of the universal
covering of a Riemann surface of finite type, 
Theorem \ref{thm:not_extremal} says the map
\begin{displaymath}
	\Teic(X) \arrow \Teic(\Delta)
\end{displaymath}
is a {\em contraction} for the Teichm\"uller metric.

\vskip 8pt

\noindent{\bf Poincar\'e series.}
To check that a map between Teichm\"uller spaces is a 
contraction,
it suffices to show the derivative of the map is a 
contracting
operator.  It is actually more convenient to work with the
{\em coderivative}, dual to the derivative operator.

Let $Q(X)$ denote the space of
holomorphic quadratic differentials $\phi (z) dz^2$ on $X$, 
such that
\[      ||\phi|| = \int_X |\phi(z)| \; |dz|^2 < \infty .
\]
With the above norm, $Q(X)$ is a Banach space.

$Q(X)$ may be naturally identified with the {\em cotangent 
space} 
to Teichm\"uller space at $X$; its norm is the infinitesimal
form of the Teichm\"uller cometric.  (As $Q(X)$ generally
fails to be a Hilbert space, this is a Finsler metric rather
than a Riemannian metric.  Technically, $Q(X)$ is the 
predual
to the tangent space.)

If $Y \arrow X$ is a covering space, then there is a 
natural 
push-forward operator
\[      \Theta : Q(Y) \arrow Q(X) .
\]
This operator is the coderivative at $X$ of
the inclusion $\Teic(X) \includesin \Teic(Y)$.

In the case of the universal covering, the operator 
$\Theta$ is
identical with the classical {\em Poincar\'e series} 
operator
\begin{displaymath}
	\Theta(\phi) \EQ \sum_{\gamma \in \Gamma}
	\gamma^*\phi,
\end{displaymath}
which takes an integrable quadratic differential $\phi$ on
the disk and converts it into an automorphic form for
the Fuchsian group $\Gamma = \pi_1(X)$ (and therefore
an element of $Q(X)$) \cite{Poincare:fuchsian}. 

A more precise formulation of Theorem 
\ref{thm:not_extremal} is

\begin{theorem}[Kra's Theta conjecture]
$||\Theta|| < 1$ for classical Poincar\'e series.
\end{theorem}

This means the inclusion of Teichm\"uller spaces is a 
contraction
even at the infinitesimal level.  

\vskip 8pt

\noindent{\bf Amenability.}
One can characterize those coverings for which contraction
is obtained in terms of the purely combinatorial notion
of {\em amenability}.

To introduce this notion, first consider the case of
a graph (1-complex) $G$.
For any set $V$ of vertices of $G$, let $\bdry V$
denote the set of vertices at distance 1 from $V$.
(A vertex at distance 1 is connected to $V$ by an
edge but does not itself lie in $V$.)
Then the {\em expansion} $\gamma$ of $G$ is given by
\begin{displaymath}
	\gamma \EQ \inf \frac{|\bdry V|}{|V|}
\end{displaymath}
over all finite sets of vertices $V$.

If the expansion is 0, $G$ is {\it amenable}; otherwise
it is nonamenable, and the boundary of any vertex set is 
comparable in size to the set itself.
For example, there is a unique infinite tree with 
degree $d$ at each vertex; it is amenable when $d=2$
(the tree is an infinite line) and nonamenable for $d>2$
(in fact the expansion constant is $d-2$).

Let $X$ be a hyperbolic Riemann surface of finite type; then
we can choose a finite graph $G \subset X$ such that
$\pi_1(G)$ surjects onto $\pi_1(X)$.  
An {\em amenable covering} $p: Y \arrow X$
is one for which $p^{-1}(G)$ is an amenable graph.
It is easy to check that the definition is independent
of the choice of $G$.
\vskip 8pt

\noindent
{\it Remark.} 
Alternatively, a covering is amenable if there exists
a linear functional (a {\em mean})
\[      m : L^\infty ( \pi_1(X) / \pi_1(Y) ) \arrow \reals,
\]
invariant under the left-action of $\pi_1(X)$,
such that $\inf(f) \le m(f) \le \sup(f)$ for all $f$.
A normal covering is amenable if the deck transformations
form an amenable {\em group};  cf. \cite{Greenleaf:book,
Pier:book}.

\begin{theorem}
\label{thm:general}
Let $Y \arrow X$ be a covering of a hyperbolic Riemann 
surface
of finite type.  Then either
\begin{enumerate}
        \item The covering is amenable, $||\Theta|| = 1$, 
and
the induced map $\mathrm{Teich}(X)\arrow 
\mathrm{Teich}(Y)$ is a global isometry
for the Teichm\"uller metric, or
        \item The covering is nonamenable, $||\Theta|| < 
1$, and
$\mathrm{Teich}(X) \arrow \mathrm{Teich}(Y)$ is contracting.
\end{enumerate}
\end{theorem}

The universal covering is easily shown to be nonamenable,
so this theorem contains the previous two.

The proof of this result is technical, but the 
relation between relaxation and nonamenability is easy
to describe.

It is related to the idea of a chain letter or pyramid game.
To join the game on a given round you must (a) pay \$1 and
(b) get two other people to join on the next round.  
Ten rounds later you leave the game, collecting 
\$1,024 from the ``descendents'' of your two new members.
In real life someone eventually loses, but on a nonamenable
graph (which is necessarily infinite), it is possible to
enrich every vertex by drawing capital at a steady rate
from infinity. 

Similarly, when a quasiconformal map is lifted to a 
cover, we can begin relaxing it on some
compact set while creating a certain amount of
additional stress near the boundary.  Using nonamenability,
the stress can be entirely dissipated to infinity.
\vskip 8pt

\noindent
{\it Bibliographical remarks.}
Theorem \ref{thm:not_extremal}
was checked in many examples by Strebel 
\cite{Strebel:theta};
Ohtake proved $||\Theta|| = 1$ for abelian coverings
\cite{Ohtake:theta}.
Other relations between amenability and function theory 
appear in \cite{Gromov:book, Brooks:amenable, 
Lyons:Sullivan}.  
\vskip 8pt

\noindent{\bf Dependence on moduli.}
Let us focus again on the case of the universal
covering $\Delta \arrow X$.
For application to iteration, it is useful to 
know {\em how much} the map $\Teic(X) \arrow \Teic(\Delta)$
contracts the Teichm\"uller metric, because
a uniformly contracting iteration has a fixed point.

\begin{theorem}
\label{thm:dep_moduli}
For the universal covering of a Riemann surface $X$ of
genus $g$ with $n$ punctures,
\begin{displaymath}
	||\Theta|| \le C(L,g,n) < 1 ,
\end{displaymath}
where $L$ is the length of the shortest geodesic on $X$.
\end{theorem}

This theorem is immediate once $||\Theta||$ is shown to be
a continuous function on the moduli space $\cM_{g,n}$,
by compactness of the set of Riemann surfaces without
short geodesics \cite{Mumford:compact}.

Moreover, $||\Theta|| \arrow 1$ as $L \arrow 0$;
there is no uniform bound on all of moduli space.
Intuitively, a Riemann surface with a short geodesic
is degenerating towards an infinite cylinder,
whose fundamental group is
$\zed$ and whose universal cover is amenable.

\section{Hyperbolic 3-manifolds}

In this section we describe Thurston's theorem on
hyperbolic 3-manifolds and how the preceding results
give a new approach to the proof.

A {\em hyperbolic} 3{\em -manifold} $N$ is a complete 
Riemannian 
3-manifold with a metric of constant curvature $-1$.
We will only consider $N$ with finitely generated 
fundamental groups.

Up to isometry, there is a unique simply connected 
hyperbolic
3-manifold, {\em hyperbolic space} $\half^3$, which 
is topologically a 3-ball.
Hyperbolic geometry tends to conformal geometry at infinity;
for example, $\half^3$ may be
compactified by the Riemann sphere $\chat$ in such a way
that hyperbolic isometries extend to conformal maps.

Since the universal cover of $N$ is isometric to $\half^3$,
we can think of $N$ as $\half^3/\Gamma$ where 
$\Gamma$ is a {\em Kleinian group},
i.e., a discrete group of hyperbolic isometries.  
There is a maximal open set $\Omega \subset \chat$ on which
$\Gamma$ acts properly discontinuously, and from this we
form the {\em Kleinian manifold} $(\half^3 \cup 
\Omega)/\Gamma$.

Thus any hyperbolic manifold $N$ is provided with a natural
Riemann surface boundary $\bdry N = \Omega/\Gamma$.

There are two important topological properties of any
hyperbolic 3-manifold.  First, $N$ is {\em irreducible}:
any 2-sphere in $N$ bounds a ball, or equivalently
$\pi_2(N) = 0$.  This is immediate from contractibility
of $\half^3$, the universal cover of $N$.

The second property is that $N$ is {\em atoroidal};
this means any incompressible torus $T^2 \subset N$
is peripheral (homotopic to an end of $N$).  
Equivalently, 
a discrete group of hyperbolic isometries isomorphic to 
$\zed \dirsum \zed$ is {\em parabolic}, i.e., it is
conjugate to a lattice of translations acting on $\chat$
by
\begin{displaymath}
	<z \mapsto z+1, \;\;\; z \mapsto z+\tau>  
\end{displaymath}
for some choice of $\tau$.
Geometrically, such a subgroup determines 
a finite volume end of $N$ (a {\em cusp}), diffeomorphic
to $T^2 \times [0,\infty)$ but rapidly narrowing.

Thurston's theorem states that in a large
category of 3-manifolds (including, e.g., all
3-manifolds with nonempty boundary), these are the
only obstructions to existence of a hyperbolic structure.

\begin{theorem}[Thurston]
\label{thm:main}
An atoroidal Haken \RM 3-manifold is hyperbolic.
\end{theorem}

A {\em Haken manifold} is one built up inductively from
3-balls by gluing along {\em incompressible} submanifolds 
of 
the boundary. 

A piece of the boundary is {\em compressible}
if there is a simple curve, which is essential in the
boundary but bounds a disk in the manifold.  When
gluing together 3-manifolds, there is danger of creating
a reducible manifold by gluing two compressible curves
together; the two disks the curve bounds join to form
an essential 2-sphere.  A Haken manifold is always 
irreducible:
by assumption, the gluing locus is incompressible, so
this danger never arises.

{\em Remark.} 
For simplicity we will suppress consideration of
the {\em parabolic locus} $P \subset \bdry M$;
in general one designates a portion of the boundary which
is to be realized as cusps.
For example, a knot complement frequently carries a 
hyperbolic structure; the boundary of a tubular
neighborhood of the knot is a torus corresponding to
a rank 2 cusp.

\vskip 8pt

\noindent{\bf Sketch of
Thurston's proof of Theorem \ref{thm:main}.}
A finite collection of disjoint balls obviously carries 
a hyperbolic structure.  Start gluing them together
along incompressible submanifolds of their boundary.
By an orbifold trick (using Andreev's theorem, see
\cite{Morgan:uniformization}), one need only deal with
the case of gluing along entire boundary components.
At the inductive step, one has a hyperbolic realization $N$
of a 3-manifold $M$ with incompressible boundary and
gluing instructions encoded by an orientation-reversing
involution 
$	\glue: \bdry M \arrow \bdry M$.
The construction is completed by the following result:

\begin{theorem}
\label{thm:gluability}
$M/\glue$ has a hyperbolic structure if and only if
the quotient is ato\-roidal.
\end{theorem}

This key result is intimately related to iteration
on Teichm\"uller space, and it this relation that
allows the results of \S\ref{sec:qc} to be brought
to bear.

Let $M$ denote a topological (or equivalently smooth)
3-manifold.  Let\break
$GF(M)$ denote the space of hyperbolic
3-manifolds $N$, which are homeomorphic to $M$.

{\em Technical remark.} 
$N$ is required to be {\em geometrically finite}, and each
$N$ is equipped with a choice of homeomorphism to $M$
up to homotopy equivalence rel boundary.  Then, as in
Teichm\"uller theory, $N_1$ and $N_2$ are regarded as
the same point in $GF(M)$ if there is an isometry 
$N_1 \arrow N_2$ in the appropriate homotopy class.

A remarkable feature of dimension 3 is
the following:

\begin{theorem}
As long as $M$ admits at least one hyperbolic realization,
there is a $1$-$1$ correspondence
between hyperbolic structures on $M$ and conformal 
structures
on $\bdry M$, i.e.,
\begin{displaymath}
	GF(M) \isom \Teic(\bdry M) .
\end{displaymath}
\RM(The isomorphism is by $N \mapsto \bdry N$.\RM)
\end{theorem}

\noindent
{\it Remarks.}
 (1) If $\bdry M$ is empty then the hyperbolic structure 
on $M$ is
unique; this is {\em Mostow rigidity}.

(2) Let $M = S \times [0,1]$ where $S$ is a surface
of genus$\geq2$; then $M$ has at least one hyperbolic
realization (consider
a Fuchsian group).  The theorem states that given any
two Riemann surfaces $X, Y \in \Teic(S)$ there is a unique
3-manifold interpolating between them; this is {\em Bers'
simultaneous uniformization theorem} \cite{Bers:simunif}.
	
(3)
The general theorem stated above was developed by Bers,
Maskit, and others and put into final form by Sullivan
\cite{Sullivan:linefield}.

To prove Theorem \ref{thm:gluability}, we must find
in $GF(M)$ the correct geometry so that the ends to
be glued together have compatible shape.  This can
be formulated as a fixed-point problem in Teichm\"uller 
space,
as follows.

First, the {\em skinning map}
\begin{displaymath}
	 \skin : \Teic(\bdry M) \arrow
		 \Teic(\overline{\bdry M})
\end{displaymath}
is defined by forming, for each $N \in GF(M)$,
the quasi-fuchsian covering spaces for each component 
of $\bdry N$ and recording the conformal
structures on the new ends that appear.  
Then, the gluing instructions determine an isometry
\begin{displaymath}
	\glue : \Teic(\overline{\bdry M})
	\arrow \Teic(\bdry M),
\end{displaymath}
and a fixed point for $\glue \compos \skin$ corresponds to
a hyperbolic manifold $N \in GF(M)$, 
whose ends fit together isometrically
under the gluing instructions.

\vskip 8pt

\noindent{\bf Properties of the skinning map.}
A detailed treatment of the skinning map would take us
deeper into the realm of Kleinian groups than we wish to
venture at the moment; here we simply
report what we would have learned upon our return.

(1) The skinning map typically contracts the 
Teichm\"uller metric ($||d\skin_X|| < 1$ at each point $X)$.
Thus if $\glue\compos\skin$
has a fixed point $X$, it can be found by iteration:
\begin{displaymath}
	(\glue\compos\skin)^n(Y) \arrow X
\end{displaymath}
for every $Y \in \Teic(\bdry M)$.
(An important exception is the case of a 3-manifold which
fibers over the circle; this can be described as 
$(M = F \times [0,1])/\glue$, where the surface $F$ is a 
fiber
and $\glue$ is the monodromy map.  In this case $\skin$ is 
an isometry and a different approach is required;
cf. \cite{Thurston:hype2}.)  

(2) In general, the skinning map does not contract {\em 
uniformly}.
This is sensible, because
$M/\glue$ might not admit a hyperbolic structure.
A potential obstruction is a cylinder $S^1\times [0,1]$
in $M$, whose ends are glued together by $\tau$ to form
a nonperipheral torus.  Then (as noted above) $M/\tau$
cannot be hyperbolic.

(3) A 3-manifold is {\em acylindrical} if every cylinder
\begin{displaymath}
(S^1 \times [0,1], S^1 \times \{0,1\}) \;\subset\; 
(M,\bdry M)
\end{displaymath}
with its ends resting on essential curves in $\bdry M$ 
can be deformed into the boundary.  
Such a manifold will remain atoroidal for
any gluing instructions, so there is no evident
obstruction to the hyperbolicity of $M/\tau$.
One expects the iteration $(\glue\compos\skin)$ to be 
robust in this case, and our first result 
(see \cite{McMullen:iter}) is

\begin{theorem}
If $M$ is acylindrical, the skinning map is {\em uniformly}
contracting: 
\begin{displaymath}
	||d\skin_X|| < C(M)  < 1.  
\end{displaymath}
Thus $M/\glue$ has a hyperbolic structure for
any choice of gluing map $\glue$.
\end{theorem}

\noindent
{\it Remark.} 
Thurston has proved the stronger assertion that the image
of the skinning map is {\em bounded} in the acylindrical 
case.

(4) We will briefly sketch how the results of the preceding
section bear on the study of $||d\skin||$.
Suppose $M$ is acylindrical.  Here is how the map
$\skin$ appears from the point of view of Riemann surfaces.
Starting with $X$ in Teich($\bdry M$), one forms countably
many copies of its universal cover $\tilde{X}$.  These
are then patched together in a complicated but canonical
pattern (depending on the topology of $M$) to form 
$\skin(X)$, the image of $X$ under the skinning map.

Now suppose $X_1$ and $X_2$ are in $\Teic(\bdry M)$, and 
let $\phi : X_1 \arrow X_2$ be an extremal quasi-conformal
map.  We can form a new quasi-conformal map 
\begin{displaymath}
	\tilde{\phi} : \skin(X_1) \arrow \skin(X_2)
\end{displaymath}
with the same dilatation as $\phi$, as follows:
lift $\phi$ to a map between corresponding copies of
the universal covers of $X_1$ and $X_2$ and complete
by continuity.

By Theorem \ref{thm:not_extremal}, $\tilde{\phi}$ can
be relaxed to a map of lesser dilatation on each
copy of the universal cover without changing its boundary
values.  Thus the relaxed maps still fit together, and
we have shown that $\skin$ contracts the Teichm\"uller
metric.  Moreover, the refined version of this relaxation
result --- Theorem \ref{thm:dep_moduli} --- yields

\begin{theorem}
$||d\skin_X|| < c(L) < 1$ where $L$ is the length of
the shortest geodesic on $X$.
\end{theorem}

This result even holds in the cylindrical case, so long
as $M$ is not an interval bundle over a surface.
The cylindrical case requires a discussion of the
nonamenability of covers of $X$ other than the universal
cover.

(5) The preceding does not yet yield uniform contraction;
rather, it reduces the problem to a study of short 
geodesics.

This is a significant simplification, however, because a
short geodesic controls the geometry of a hyperbolic
manifold over a large distance, by the Margulis lemma.
When combined with the theory of geometric limits of
quadratic differentials 
(cf. Appendix to \cite{McMullen:amenable}),
one finds a qualitative picture which again forces
contraction {\em unless} the short geodesic lies
on one end of a compressing cylinder in the 3-manifold.

This cannot occur in the acylindrical case, so then the
contraction is uniform.  

In general we find:

\begin{theorem}
For any initial guess $Y \in \Teic(\bdry M)$, either
\begin{displaymath}
	Y_n \EQ (\glue \compos \skin)^n(Y) \arrow 
	\mbox{a (unique) fixed point $X$,} 
\end{displaymath}
or $Y_n$ develops short geodesics, bounding cylinders 
in $M$ linked by $\tau$ to form a nonperipheral torus
in $M/\glue$.  

Thus the gluing problem is solvable if and
only if $M/\glue$ is atoroidal.
\end{theorem}

\section{Epilogue}

Two other iterations on Teichm\"uller space have
been analyzed with a similar paradigm, and indeed one
of our motivations was to give a parallel treatment
of the skinning map.  
One iteration is the action on Teich($S$) of
an element of the mapping class group of the surface $S$
\cite{Thurston:mcg, FLP}; 
the other arises
in the construction of critically finite rational maps
\cite{Thurston:rational, Douady:Hubbard:Thurston}.

The paradigm is to play (a) geometric control
in the absence of short\break
geodesics, against (b) topological
information in their presence.  

For example, a mapping class whose minimum translation
distance is achiev-\break
ed is {\em geometric}; it either fixes
a point in Teichm\"uller space (and can be represented by
an {\em isometry}), or it fixes a geodesic (this is the 
{\em pseudo-Anosov} case.)  Otherwise the minimum 
translation
distance is not achieved, so there are Riemann surfaces
tending to infinity in moduli space and translated
a bounded distance.  In the end of moduli space the surfaces
have short geodesics, which are necessarily permuted; 
thus we have the topological conclusion that the mapping
class is {\em reducible}.  This approach appears in
\cite{Bers:mcg}.

Similarly, given a smooth branched cover $f : S^2 \arrow 
S^2$
whose critical points eventually cycle, one seeks a 
conformal 
structure preserved up to isotopy by $f$ (rel the 
post-critical 
set); then $f$ can be geometrized as a rational map.
This desired conformal structure can again be described as
a fixed point for an iteration on Teichm\"uller space.
The iteration contracts uniformly in the absence of short
geodesics; in their presence,
one locates a topological obstruction to geometrization.

Finally, in retrospect it is promising to think
of the skinning map as an instance of {\em renormalization};
a similar approach to the Feigenbaum phenomenon
(using infinite-dimensional Teichm\"uller spaces) has
been proposed by Douady and Hubbard and 
pursued by Sullivan \cite{Sullivan:icm}.

\vskip 8pt

\noindent{\bf Addendum October 1991.}
The connection with renormalization is now
understood more precisely; construction of
the Feigenbaum fixed-point closely resembles the
construction of hyperbolic structures on
3-manifolds, which fiber over the circle---the
one case omitted from the discussion above.
This analogy and others (which relate the construction
of rational maps to that 
of hyperbolic 3-manifolds) are
discussed in more detail in
\cite{McMullen:icm90}.
Abundant progress towards understanding 
renormalization of quadratic polynomials appears
in \cite{Sullivan:renormalization}.

A new approach to the Theta conjecture and its
variants is discussed in \cite{Barrett:Diller,
McMullen:newamen}.



\end{document}